\documentclass{article}
\usepackage{epsf}
\usepackage{times}
\usepackage{graphicx}

\DeclareSymbolFont{AMSb}{U}{msb}{m}{n}
\DeclareSymbolFontAlphabet{\Bbb}{AMSb}
\DeclareMathSymbol{\bba}{\mathbin}{AMSb}{'101}
\DeclareMathSymbol{\bbr}{\mathbin}{AMSb}{'122}
\DeclareMathSymbol{\bbc}{\mathbin}{AMSb}{'103}
\DeclareMathSymbol{\bbh}{\mathbin}{AMSb}{'110}
\DeclareMathSymbol{\bbo}{\mathbin}{AMSb}{'117}

\begin{document}

\title{
A Very Intuitive Geometric Picture of the 24-cell, $E_8$ and $\Lambda_{16}$ Lattices Given by Using the Hopf Maps
}

\author{
 {\bf Eric Lewin Altschuler$^1$ and Antonio P\'erez--Garrido$^2$}\\
$^1$Department of Physical Medicine \& Rehabilitation, UMDNJ\\
\small University Hospital, 150 Bergen Street, B-403, Newark, NJ 07103, USA\\
\small email: eric.altschuler@umdnj.edu\\
$^2$Departamento de F\'\i sica Aplicada, UPCT\\
\small Campus Muralla del Mar, Cartagena, 30202 Murcia, Spain\\
\small email:Antonio.Perez@upct.es
}
\maketitle

\begin{abstract}
 We use the Hopf fibrillation to give simple and intuitive geometric constructions of the 24-cell, $E_8$ and 
 $\Lambda_{16}$ lattices.
 \end{abstract}
 
 The $n$-dimensional {\it kissing} problem asks how many non overlapping unit $n$-dimensional spheres can be placed touching a central unit 
 n-dimensional sphere?  This question has applications to making efficient codes and other problems
 \cite{CS99}. The two dimensional sphere is the disk of points $x^2 + y^2=1$, and by 
 using $N$ identical circular coins the reader can easily convince oneself that the kissing number in two dimensions ($K_2$) 
 is  six.   Similarly, the one dimensional sphere is a line segment, and $K_1 = 2$.  Already in three dimensions the problem
  becomes   much more interesting.  Indeed, the Seventeenth century featured a dispute between Isaac Newton who 
  believed 
  that $K_3 =   12$, and David Gregory who thought that $K_3 = 13$.
   Perhaps not surprisingly, Newton was correct but it took more than two and one half centuries to prove this\cite{SW53}.  
   Using clever  linear programming arguments it was proven some decades ago that $K_8 = 240$ and $K_{24} = 196560$
   \cite{Le79,OS79}    
   and that $K_4 = 24$ or $25$.  Recently, a proof that $K_{24} = 24$ using much more extensive and subtle use of 
  linear programming has be presented\cite{Mu03}, 
  however, the proof has not yet been fully vetted \cite{PZ04}.
  There exists a map known as the (first) {\it Hopf} map between the 
   surface of a    three dimensional sphere and the surface of a four dimensional sphere, which essentially constructs the 
   surface of a four 
  dimensional sphere by considering there to be a circle of points at every point on the surface of a three dimensional 
  sphere.   Here we emphasize that the Hopf map  also gives a very intuitive 
way of appreciating the $N = 24$ kissing configuration for $S^3$--known as the {\it 24-cell} and then use the second and third Hopf maps to give intuitive descriptions of the $E_8/K_8$ lattice, and the so called $\Lambda_{16}$ 
lattice, currently the best known kissing configuration in 16 dimensions.

The surface of a four dimensional sphere (a three dimensional locus or manifold also known as $S^3$) is defined as the points $x^2+y^2 + z^2+w^2 = 1$ (as the surface of a three dimensional sphere ($S^2$) is defined as the points 
$x^2+y^2 + z^2 = 1 $).  The Hopf map is one between the points on the surface of a four dimensional sphere and the pair of complex numbers $(w, z)$ with $|w|^2 + |z|^2 = 1$ 
\begin{equation}
(w,z) \rightarrow (2 w z^*, |z|^2-|w|^2)\,  {\rm in}\,\,  \bbc\times\bbr = \bbr^3.
\end{equation}
One easily checks that:
\begin{equation}
|2 w z^*|^2 + (|z|^2-|w|^2)^2
= 4 |w|^2 |z|^2 + (|z|^2-|w|^2)^2 = (|z|^2+|w|^2)^2 = 1.
\end{equation}
So this does map to the ordinary sphere. If one fixes a point of the
ordinary sphere say  $(a,t)$  where  $a$  is complex, $t$ is real and
$|a|^2 + t^2 = 1$, then its {\it fiber}, i.e., the set of all points
which map to it, is a circle
\begin{equation}
\left(\frac{a e^{i \theta}}{\sqrt{2(1+t)}}, e^{i \theta}\sqrt{(1+t)/2}\right).
\end{equation}
Further details, discussions and proof 
of the Hopf map from $S^2$ to $S^3$ is given in\cite{Ho31}.

The only known configuration on $S^3$ with twenty-four kissing spheres is the so-called 24-cell. In this convex four dimensional polytope all of the faces are  octahedra.  The standard way of representing the 24-cell is by the coordinates: 
   \[
\left(\pm \frac{\sqrt{2}}{2}, \pm \frac{\sqrt{2}}{2} ,0,0 \right),\,\,
\left(\pm \frac{\sqrt{2}}{2},0, \pm \frac{\sqrt{2}}{2} ,0 \right),\,\,
\left(\pm \frac{\sqrt{2}}{2},0,0,  \pm \frac{\sqrt{2}}{2} \right)
\]

\[
\left(0, \pm \frac{\sqrt{2}}{2}, \pm \frac{\sqrt{2}}{2} ,0 \right),\,\,
\left(0, \pm \frac{\sqrt{2}}{2},0, \pm \frac{\sqrt{2}}{2}  \right),\,\,
\left(0, 0, \pm \frac{\sqrt{2}}{2}, \pm \frac{\sqrt{2}}{2}  \right)
\]
 
 This configuration is also known as the $E_4$ lattice, and one can form an analogous lattice in any dimension.  By rotating 
 on $S^3$ one can also easily see that these same twenty-four centers of the kissing spheres can be obtained by lifting, via 
 the (first) Hopf map from six points on $S^2$ arranged with one at each pole and four arranged at ninety degree angles 
 around the equator:

Place one point on each pole of a three dimensional sphere and four equally spaced points on the 
 equator, i.e. the vertex of an octahedron. These points are at the {\it antipodal} points of the three axes  of $S^2$. 
 These points can be expressed as $(a,t)$, as stated above:
 \[
 (0,1),\, (0,-1),\, (1,0),\, (-1,0),\, (i,0),\, (-i,0) 
 \]
 
Then four 
 points can be placed on each circle 
via the Hopf map to give a total of twenty four points. 
\[
\left(e^{i\theta},0\right), \,\, \theta=\frac{\pi}{4}(2k+1), \, k=0,1,2,3
\] 

\[
\left(0,e^{i\theta}\right), \,\, \theta=\frac{\pi}{4}(2k+1), \, k=0,1,2,3
\] 

\[
\left( \frac{e^{i\theta}}{\sqrt{2}}   , \frac{e^{i\theta}}{\sqrt{2}}          \right), \,\, \theta=\frac{\pi}{2}k, \, k=0,1,2,3
\]

\[
 \left(   -\frac{e^{i\theta}}{\sqrt{2}}   , \frac{e^{i\theta}}{\sqrt{2}}          \right), \,\, \theta=\frac{\pi}{2}k, \, k=0,1,2,3
\] 

\[
 \left(   i\frac{e^{i\theta}}{\sqrt{2}}   , \frac{e^{i\theta}}{\sqrt{2}} \right) \equiv  
  \left(   \frac{e^{i(\theta+\pi/2)}}{\sqrt{2}}   , \frac{e^{i\theta}}{\sqrt{2}}       \right), \,\, \theta=\frac{\pi}{2}k, \, k=0,1,2,3
\] 

 It is easy to see that these twenty-four points on $S^3$ (the 
surface of a unit four-dimensional sphere) are the same as the points of a 24-cell.

\begin{figure}
\begin{center}
\leavevmode
\includegraphics[angle=0,width=9.5cm]{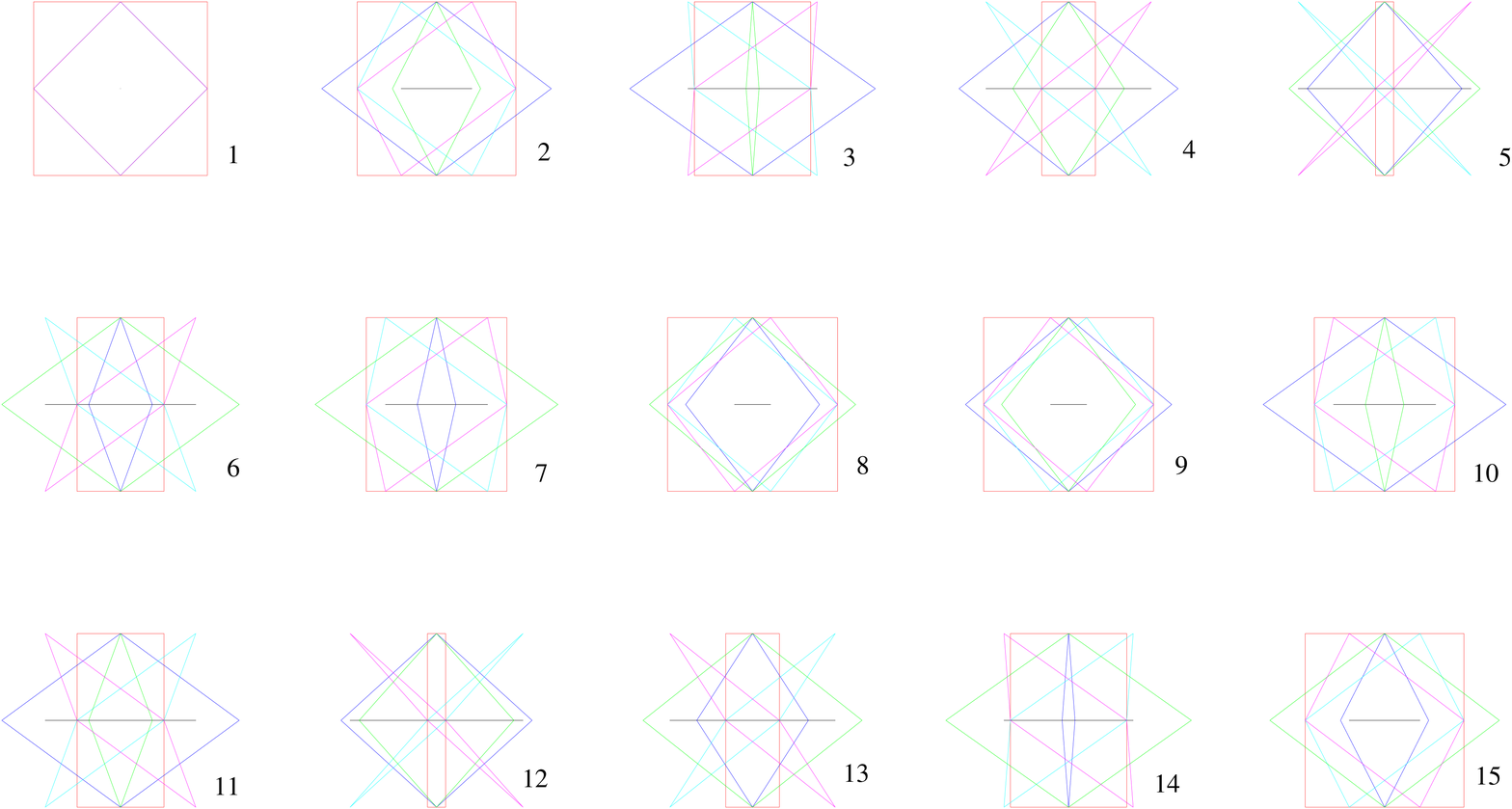}
\end{center}
\caption{24-cell. Parallel projection to a 2D plane. 
Several views during an $i\cdot2\pi/15$  rotation in 4D  are plotted. Each color corresponds to a different circle on $S^2$. 
The symmetry of the configuration is evident.  The antipodal construction is illustrated as the circles (in $S^2$) 
are ninety or one hundred eighty degrees from each other.  
That nearest neighbors are points on different circles
(see text for discussion) is illustrated, for example,
in  subfigure 1 where
 one sees that the black circle, actually seen as a single point, is surrounded by 
 all other circles but the red circle which is its antipodal. 
}
\end{figure}

\begin{figure}
\begin{center}
\leavevmode
\includegraphics[angle=-90,width=9.5cm]{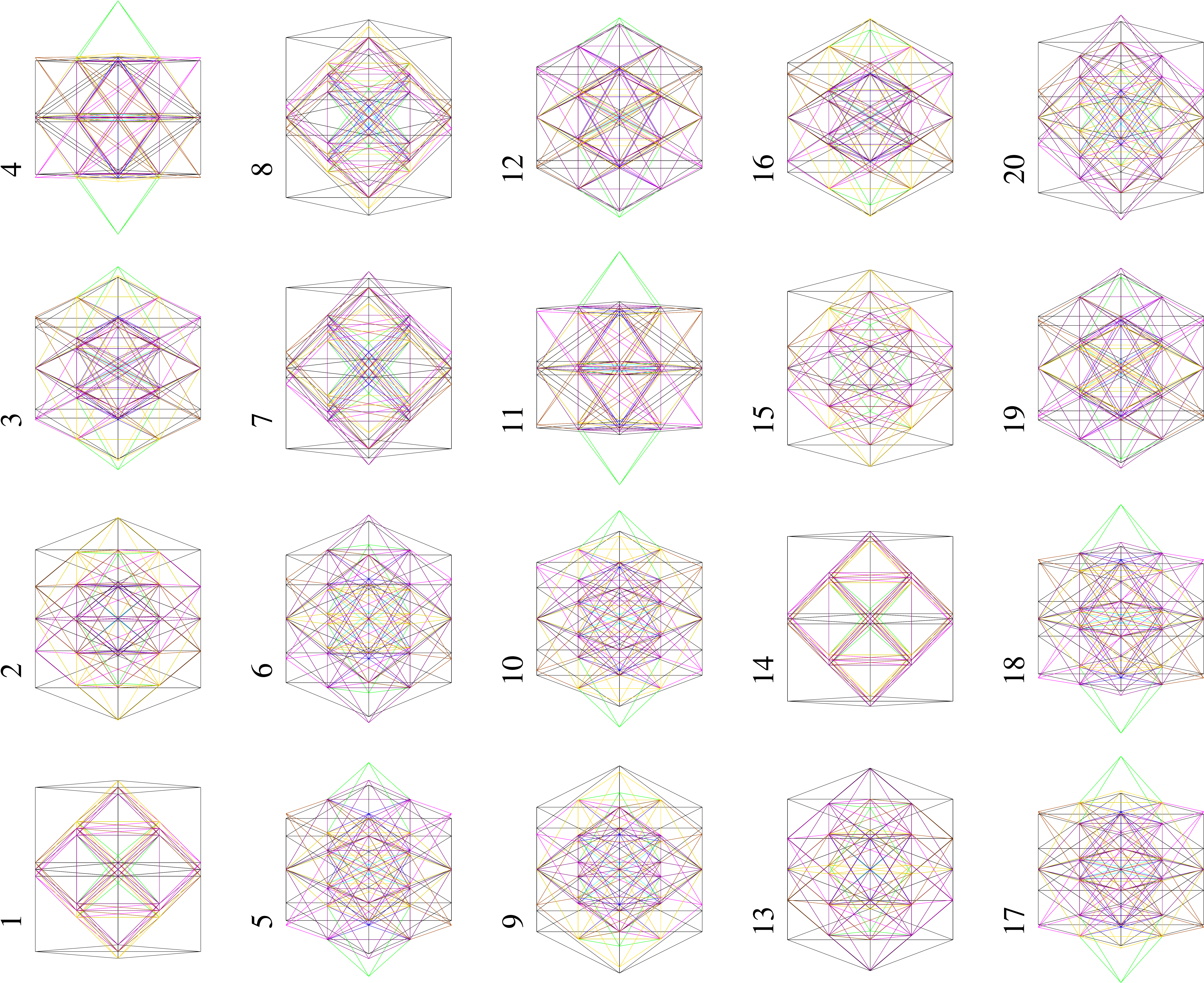}
\end{center}
\caption{
$E_8$. Parallel projection to a 2D plane of several views 
during a rotation by $i\cdot2\pi/20$  around a plane in 8D are plotted.
Each color corresponds to a different circle on $S^4$.
Again we see the symmetry inherent in this configuration, and the 
  antipodal construction being highlighted as the circles ($S^3$) are ninety or 
180 degrees apart. 
}
\end{figure}

The twenty-four points can 
trivially be seen to be a distance greater or equal than 
one from each other and thus a kissing configuration in four dimensions of twenty-four 
points. This is the known configuration  illustrating that $K_4$ is $\geq 24$ and is
also know as the {\it 24-cell} because those points are the vertices of a 4D polytope made of 24 octahedra.  
4D polytopes are assemblies of 3D polyhedra (cells) as 3D polyhedra are assemblies of 2D polygons.

Also, the Hopf map gives an extremely simple and somewhat intuitive way to 
think about the 24-cell: For one point on each equator and four spaced at ninety 
degrees along the equator then each point in $S^2$ is a distance of ninety 
degrees away from its nearest neighbor.  Now associate with each of these six points 
a {\it Hopf circle} placing four points again spaced ninety degrees.  The 
symmetry of the configuration becomes self-evident.  

Points on the same circle are ninety degrees--or distance $\sqrt{2}$ apart from 
each other.  Each point has eight nearest neighbors--two points each from the 
four non-antipodal circles--that are sixty degrees or distance one apart.  

Now, in addition to the first Hopf map from $S^3$ to $S^2$, there is a second Hopf map from $S^7$ to $S^4$
  (and a third Hopf map from  $S^{15}$ to $S^8$).  Intuitively the second (third) Hopf map uses quaternions (octonions) to 
  accomplish the map.  
  Using the  Cayley-Dickson construction for the normed division algebras, $\bba_0=\bbr$ (real)   $\bba_1=\bbc$ (complex) 
  $\bba_2=\bbh$ (quaternions)  and $\bba_3=\bbo$ (octonions), we can make $\bba_n$ from $\bba_{n-1}$ for $n=1,2,3$.
  In this construction, an element in $\bba_n$ is made of a pair of elements $a,b\in\bba_{n-1}$ with the multiplication:
  \begin{equation}
  (a,b)(c,d)=(ac-db^*,a^*d+cb)
  \end{equation}
  and the conjugation in $\bba_n$:
  \begin{equation}
  (a,b)^*=(a^*,-b)
  \end{equation}
It is also possible to continue this procedure to get $\bba_4$ (sedenions) but it is no longer a division algebra.
Hopf maps can be compactly defined as a map $h_1$ from $\bba_n\otimes\bba_n$ to $\bba_n\cup\lbrace\infty\rbrace$
followed by a second map $h_2$ from $\bba_n\cup\lbrace\infty\rbrace$ to $S^{2^n}$ (stereographic projection), for $n=0,1,2,3$\cite{MD01}. $h_1$ and $h_2$ can be stated as:
\begin{equation}
h_1:\,\,(a,b) \,\, \longrightarrow \,c=ab^{-1}  
\end{equation}
where $a,b\in\bba_n$ and $\left| a \right|^2+\left| b \right|^2=1$, so the point $(a,b)$ is actually a point of $S^m$, being $m=2^{n+1}-1$.  
\begin{equation}
 h_2:\,\,c \,\, \longrightarrow \,X_i\,\, (i=1,\ldots,2^n+1), \,\,\, \sum_{i=1}^{2^n+1}X_i^2=1
 \end{equation}
 
Now, Dixon (see \cite{Di95} and refs. therein) seems to have been the first to have appreciated not only that the first Hopf map 
can be used to generate the 24-cell from points on $S^2$, but that $E_8$ is generated from ten 24-cell's lifted from $S^4$, and $\Lambda_{16}$ is generated from 
18 $E_8$ lattices lifted from $S^8$.  Dixon's discussion was mostly topological.  Here we give a remarkably simple and intuitive geometric construction for $E_8$ and $\Lambda_{16}$.  
We  note (see \cite{MD01} and refs. therein) that such Hopf make constructions may have utility in quantum computers and communication.

$S^7$ is the surface of an eight dimensional sphere.  The kissing configuration of eight dimensional spheres on the surface of an eight dimensional sphere as mentioned is known to be 240 points arranged in the $E_8$ lattice. The kissing configuration in five dimensions ($K_5$) which is points arranged on the surface of $S^4$ is thought to be 40 points arranged in an $E_5$ lattice, but there is no proof of this.  Initially we wondered if we could use the second Hopf map to lift from the forty kissing points on $E_5$ six points each onto $S^7$ and obtain the $E_8$ lattice/$K_8$ configuration.  We have not been able to do this, however, we noticed that by taking again the 10 {\it antipodal} points from the axes on $S^4$,  $(\pm1,0,0,0,0),$,    $(0,\pm1,0,0,0),$, $(0,0,\pm1,0,0),$
$(0,0,0\pm1,0)$
and  $(0,0,0,0,\pm1),$ and lifting 10 24-cells to $S^7$ we get the $E_8$ lattice. 
   This construction is illustrated in Figure 2.  Again, our construction 
immediately illustrates that as for the 24-cell points are separated by sixty, 
ninety, 120 or 180 degrees.  
  Our construction also intuitive explains why each point has 56 nearest 
neighbors: eight on "its own" $S^3$ circle, and six (two per orthogonal axis of a 3-dimensional object) on each of the eight other non-antipodal circles on $S^4$.

Similarly, by lifting from the 16 antipodal points of the axes of $S^8$ an $E_8$ lattice one gets the $\Lambda_{16}$ lattice! 
This construction again illustrates why $\Lambda_{16}$ is 
a kissing configuration with points separated by angles sixty, ninety, 120 or 
180 degrees, and why each point on $\Lambda_{16}$ has  280 nearest neighbors = 56 on 
"its own" $S^7$ circle + 14 (two points per orthogonal axis of a 7-dimensional object) $\times$ 16,  
where 16 is the number of non-antipodal $S^7$ circles on $S^8$.
This construction suggests, but of course in no way proves, that $\Lambda_{16}$ may be a configuration of maximum kissing number. Our construction may give a 
helpful picture in proving, or disproving this, or in attacking other problems such as finding the maximal packing configurations in 4, 8 or 16 dimensions.

We thank Andrew Gleason and Richard Stong for helpful discussions. APG was partly supported by 
Spanish MCyT under grant No. MAT2003--04887.

\bibliographystyle{nature}

\end{document}